\documentclass[11pt]{article}
\usepackage{amsfonts, amsmath, amssymb, amsthm}
\usepackage[T1,T2A]{fontenc}
\usepackage[utf8]{inputenc}

\PassOptionsToPackage{obeyspaces}{url}
\usepackage[colorlinks=true,citecolor=blue,urlcolor=blue,linkcolor=blue,bookmarksopen=true]{hyperref}
\usepackage{breakurl}
\usepackage{tikz}

\usepackage{etoolbox} 
\patchcmd{\thebibliography}{\leftmargin\labelwidth}{\leftmargin\labelwidth\addtolength\itemsep{-0.1\baselineskip}}{}{}

\author{Boris Bukh\thanks{Department of Mathematical Sciences, Carnegie Mellon University, Pittsburgh, PA 15213, USA. Supported in part by Sloan Research Fellowship and by U.S.\ taxpayers through NSF CAREER grant DMS-1555149. Part of the work was done during the visit
to Universit\'e Paris-Est Marne-la-Vall\'ee supported by LabEx B\'ezout (ANR-10-LABX-58).} \and Alfredo Hubard \thanks{Laboratoire d'Informatique Gaspard-Monge, Universit\'e Paris-Est Marne-la-Vall\'ee, France. alfredo.hubard@u-pem.fr.}}

\title{On a topological version of Pach's overlap theorem}
\date{}

\newtheorem{theorem}{Theorem}
\newtheorem{lemma}[theorem]{Lemma}

\newtheorem{definition}[theorem]{Definition}

\theoremstyle{remark}

\newcommand*{\eqdef}{\stackrel{\text{\tiny{def}}}{=}}            
\newcommand*{\abs}[1]{\lvert #1\rvert}                           
\newcommand*{\norm}[1]{\lVert #1\rVert}                          
\newcommand*{\veps}{\varepsilon}                                 
\newcommand*{\R}{\mathbb{R}}                                     
\newcommand*{\RP}{\mathbb{RP}}                                   
\newcommand*{\Sph}[1][d]{\mathbb{S}^{#1}}                        
\newcommand*{\F}{\mathbb{F}}                                     
\newcommand*{\M}{M}                                              

\def\D#1-{$\if\relax\detokenize{#1}\relax d\else (d#1)\fi$\nobreakdash-\hspace{0pt}}

\DeclareMathOperator{\sign}{sign}                                
\DeclareMathOperator{\supp}{supp}                                
\DeclareMathOperator{\prj}{prj}                                  

\newcommand*{\TODO}[1]{}

\begin{document}
\maketitle

\begin{abstract}
Pach showed that every $d+1$ sets of points $Q_1,\dotsc,Q_{d+1} \subset \R^d$ contain linearly-sized subsets $P_i\subset Q_i$ such that
all the transversal simplices that they span intersect. We show, by means of an example, that a topological 
extension of Pach's theorem does not hold with subsets of size $C(\log n)^{1/(d-1)}$. We show
that this is tight in dimension $2$, for all surfaces other than $\Sph[2]$. Surprisingly, the optimal
bound for $\Sph[2]$ in the topological version of Pach's theorem is of the order $(\log n)^{1/2}$. We conjecture that, 
among higher-dimensional manifolds, spheres are similarly distinguished. This improves upon the results of B\'ar\'any, Meshulam, Nevo and~Tancer.
\end{abstract}

\section*{Introduction}
\paragraph{Pach's overlap theorem.}
Let $Q_1,\dotsc,Q_{d+1}$ be $d+1$ sets of points in~$\R^d$. A family of subsets $P_1,\dotsc,P_{d+1}$
with $P_i\subset Q_i$ is a \emph{Pach family} if all the simplices of the form $p_1\dotsc p_{d+1}$ with $p_i\in P_i$
have a common intersection point. Pach's overlap theorem \cite{pach}
asserts that if each set $Q_i$ is of size $n$ then there is a Pach family with $\abs{P_i}\geq c_d n$. The constant $c_d$ depends only on the dimension and its optimal value is unknown. The same-type lemma of B\'ar\'any and Valtr \cite{barany_valtr} is closely related to Pach's result. 

\paragraph{Topological version.}
Many results in combinatorial geometry admit a topological generalization by replacing linear functions with
arbitrary continuous functions. A celebrated example concerns planar graphs: a graph can be drawn in the plane
with straight edges if and only if it can be drawn with `curvy' edges. Closer to the subject of the
present paper are the topological Tverberg theorem of B\'ar\'any, Shlosman and  Sz{\H{u}}cs \cite{topological_tverberg}
and Gromov's topological overlap theorem \cite{gromov}, which generalize Tverberg's theorem \cite{tverberg} and
B\'ar\'any's overlap theorem \cite[Theorem~5.1]{barany_caratheodory} respectively.

It is natural to ask if Pach's overlap theorem admits such a generalization. In order to state it, we need some standard notation.
For a simplicial complex $Y$, let $Y^{\leq d}$ be the $d$-skeleton of $Y$, i.e., the complex consisting of all the simplices of dimension~$\leq d$.
If $Y_1,Y_2$ are simplicial complexes, $Y_1\ast Y_2\eqdef\{A\sqcup B: A\in Y_1,B\in Y_2\}$ denotes their join (where $\sqcup$ denotes the disjoint union). In the following, we let $V_1,V_2,\dotsc,V_{d+1}$ be $d+1$ point sets of size $n$ each. By a slight abuse of notation, we treat each $V_i$ also as a
$0$-dimensional simplicial complex. 
Put \[X\eqdef V_1\ast\dotsb\ast V_{d+1}.\]
The complex $X$ consists of $d$-simplices, which contain one vertex from each $V_i$, and their faces.

\begin{definition}
Consider a continuous map $\phi\colon X\to \M$ into a $d$\nobreakdash-dimensional manifold. A \emph{topological Pach family for $\phi$} consists of sets $P_i\subset V_i$
and a point $p\in M$ such that the $\phi$\nobreakdash-image of every $d$-simplex in $P_1\ast \dotsb \ast P_{d+1}$ contains~$p$. 
Let $f_{\M}(n)$ be the largest integer such that, for every continuous map $\phi\colon X\to \M$, there is 
a topological Pach family all of whose sets are of size~$f_{\M}(n)$. 
\end{definition}

When $\M=\R^d$ and $\phi$ is an affine map, this notion reduces to the usual (non-topological) Pach family. In this paper, the map $\phi$ is an arbitrary continuous map
into a topological manifold $\M$ which is assumed to have no boundary, but is not necessarily compact.
In \cite{barany_meshulam_nevo_tancer}, it was shown that, unlike in the affine case, a topological Pach family of linear size need not exist. 
Namely, \cite{barany_meshulam_nevo_tancer} gave a construction demonstrating that $f_{\R^d}(n)\leq C_d n^{1/d}$;
the same paper also showed that $f_{\R^d}(n)\geq c_d (\log n)^{1/d}$. We present a better construction and improve the lower bound in dimension~$2$, showing
that our construction is tight in this case.

Below and everywhere else in the paper the logarithms are to base $e=2.71\dots$.
\begin{theorem}\label{thm:top}
For every manifold $\M$ of dimension~$d\geq 2$,  
\begin{align*}
f_{\M}(n)&\leq 30 (\log n)^{1/(d-1)},\\
\intertext{for $n$ large enough. This is sharp for $d=2$: for every surface $\M\neq \Sph[2]$,}
f_{\M}(n)&\geq 10^{-14} \log n.
\end{align*}
\end{theorem}

Much to our surprise the asymptotics of $f_{M}(n)$ is sensitive to the target space.
The asymptotics of $f_{\Sph[2]}(n)$ is different from the asymptotics of $f_M(n)$ for any other surface $\M\neq \Sph[2]$!

\begin{theorem}\label{thm:sphere}
For each $d\geq 1$, \[\alpha_d \leq \frac{f_{\Sph}(n)}{(\log n)^{1/d}}\leq 2,\]
for $n$ large enough. Here, $\alpha_d= \frac{2^{-d^2-1}}{(d+1)!}$.
\end{theorem}

We begin with the proof of the easiest result, which is Theorem~\ref{thm:sphere}. We then explain the similar, but more complicated
Theorem~\ref{thm:top}. We finish this paper with some remarks and open problems.

\section*{The case $\M=\Sph$: proof of Theorem~\ref{thm:sphere}}
\paragraph{Construction.} We identify $\Sph$ with the one-point compactification of $\R^d$.
Map the vertices of $X$ to a set in general position, and then extend to the $(d-1)$-skeleton by linearity.
For each $d$-face of $X$, there are two natural ways of extending $\phi$ to that face. We may stay inside $\R^d$ and 
fill in the simplex by extending $\phi$ linearly, or we may use the point at infinity to obtain the `inverted' filling. We choose one of these two ways
uniformly at random, for each $d$-face.

\begin{center}
\begin{tikzpicture}
\begin{scope}
\filldraw[draw=black, fill=gray, line width=1mm] (-0.5,-0.288675) -- (0.5,-0.288675) -- (0,0.57735) -- cycle;
\end{scope}
\begin{scope}[xshift=4cm]
\clip plot[no marks] file {polygon.txt};
\filldraw[even odd rule, draw=black, fill=gray, line width=1mm] (-4,-4) rectangle (4,4)  (-0.5,-0.288675) -- (0.5,-0.288675) -- (0,0.57735) -- cycle;
\end{scope}
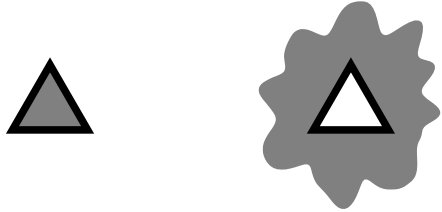
\node[align=left] at (2,-1.7) {\scriptsize The usual (left) and the `inverted' (right) fillings of a triangle};
\end{tikzpicture}
\end{center}

The $(d+1)n$ vertices of $X$ span
$\binom{(d+1)n}{d}$ hyperplanes, which in turn partition $\Sph$ into $O(n^{d^2})$~regions. The points in each such 
region are covered by the same set of $\phi$-images of $d$-faces of~$X$. 
Hence,
the common intersection point of every Pach family can be taken to be a center of one of the regions.
Let $m\eqdef 2(\log n)^{1/d}$.  Given a point $p\in \Sph$ and a \D+1-tuple of $m$-element sets $(P_1,\dotsc,P_{d+1})$,
the probability that $p$ is a common intersection point of all the simplices spanned by $P_1,\dotsc,P_{d+1}$ is
$2^{-m^{d+1}}$ since the fillings are chosen independently. By the union bound over all tuples $(P_1,\dotsc,P_{d+1})$ and
over all region centers $p$, it follows that a Pach family of size $m$ exists with probability at most
$O\Bigl(\binom{n}{m}^{d+1} n^{d^2} 2^{-m^{d+1}}\Bigr)=O\Bigl((n^{d+1}2^{-m^d})^mn^{d^2}\bigr)=O(n^{-m}n^{d^2})=o(1)$. In particular, for some choice of fillings, a Pach family of this
size does not exist.


\paragraph{Lower bound.} This argument is from \cite{barany_meshulam_nevo_tancer}, so we only sketch it.
Denote by $X^{=d}$ the set of \D-simplices in~$X$. By Gromov's topological overlap theorem \cite{gromov}, for any continuous map $\phi\colon X\to \Sph$,
there is a point $p$ that is in the $\phi$-images of $\alpha_d \abs{X^{=d}}$ many $d$-faces of~$X$, where $\alpha_d>0$
is an absolute constant. Let $H=\{\sigma\in X^{=d} : p\in \phi(\sigma)\}$. Since each $d$-simplex contains exactly one vertex from each~$V_i$,
we can regard $H$ as \D+1-partite \D+1-uniform hypergraph with parts~$V_1,\dotsc,V_{d+1}$.
Every \D+1-partite \D+1-uniform hypergraph of density $\alpha>0$ contains a complete \D+1-subhypergraph all of whose parts are
of size $\lfloor \alpha (\log n)^{1/d}\rfloor$. This is Nikiforov's explicit estimate \cite{nikiforov_kst} for Erd\H{o}s's generalization \cite{erdos}
of the K\"ovari--Sos--Tur\'an theorem. 

The explicit value for $\alpha_d$ of $1/(d+1)!2^{d^2+1}$ will be obtained in the next section, after a more careful discussion
of Gromov's theorem. Plugging this value into Nikiforov's estimate, we obtain the desired lower bound.

\section*{The case $\M\neq \Sph$: proof of Theorem~\ref{thm:top}}
\paragraph{Construction.} This is similar to the construction for $\Sph$, except that we
cannot randomize the \D-faces of~$X$. We randomize the \D{-1}-faces instead. 

We say that $H$ is a \emph{$(d,n)$-graph} if $H$ is a complete
\D+1-partite $d$-uniform hypergraph in which each part is of size~$n$.
In our probabilistic construction we will require an upper bound on the probability that a random subgraph of a $(d,n)$-graph contains
no \D+1-clique. It is possible to prove a very precise bound 
using hypergraph containers (see \cite[Corollary 2.4]{saxton_thomason} or \cite[Theorem 8.1]{balogh_morris_samotij} for similar non-partite
results). Instead of doing that, we opt for a simple proof of a cruder upper bound using a restatement of the Loomis--Whitney inequality
\cite{loomis_whitney} due to Bollob\'as--Thomason\cite{bollobas_thomason}.

For a set family $\mathcal{F}\subset 2^X$ and a set $Y\subset X$, we define the \emph{trace on $Y$} by
$\mathcal{F}|_Y\eqdef \{S\cap Y : S\in \mathcal{F}\}$.
\begin{lemma}[Theorem~7 from \cite{bollobas_thomason}]\label{lem:loomis_whitney}
Let $E$ be a finite set, and $\mathcal{C}\subset 2^E$ be a family of subsets of $E$ such that each $x\in E$ is
in equally many sets of~$\mathcal{C}$. Let $F$ be a uniformly random subset of~$E$. 
Suppose that a family $\mathcal{F}\subset 2^E$ and a constant $c<1$ satisfy
$\Pr[F\cap Y \in \mathcal{F}|_Y]\leq c^{\abs{Y}}$ for every $Y\in \mathcal{C}$. Then
$\Pr[F\in \mathcal{F}]\leq c^{\abs{E}}$.
\end{lemma} 
\begin{lemma}
Each $(d,n)$-graph admits a coloring of its edges into two colors
so that each $(d,m)$\nobreakdash-subgraph contains monochromatic \D+1-cliques
in each of the two colors, for all $m\geq 25 (\log n)^{1/(d-1)}$. 
\end{lemma}
\begin{proof}
Let $H$ be any $(d,m)$-subgraph of the $(d,n)$-graph. Let $E$ be the edge set of~$H$, let $\mathcal{F}$ be the set of
all subgraphs of $H$ that contain no $(d+1)$-clique, and let $\mathcal{C}$ be the set of all $(d+1)$-cliques in $H$.
Then $\Pr[F\cap Y\in\mathcal{F}|_Y]=1-2^{-d-1}$ for a random $F\subset E$ and any $(d+1)$-clique~$Y$. Because each edge
is covered by equally many $(d+1)$-cliques and $\abs{E}=(d+1)m^d$, Lemma~\ref{lem:loomis_whitney} implies that $\Pr[F\in \mathcal{F}]\leq 
(1-2^{-d-1})^{m^d}$. Hence the probability
that $H$ does not satisfy the conclusion in a random coloring of the $(d,n)$\nobreakdash-graph is
at most $2(1-2^{-d-1})^{m^d}$. 
Since $\binom{n}{m}^{d+1}\cdot 2(1-2^{-d-1})^{m^d}<n^{m(d+1)}\cdot 2\exp(-2^{d+1}m^d)<1$ for $m=\lceil 25(\log n)^{1/(d-1)}\rceil$,
by the union bound over all $(d,m)$-subgraphs $H$ it follows that the probability of some $H$ failing the conclusion is $<1$, and so
the desired coloring exists.
\end{proof}

Recall that $X=V_1\ast\dotsb\ast V_{d+1}$. The $(d-1)$-simplices of $X$ form a $(d,n)$\nobreakdash-graph. Choose a $2$\nobreakdash-edge-coloring $\chi\colon X^{=(d-1)}\to \{1,-1\}$ of this hypergraph
as in the preceding lemma.\smallskip

\parshape=13 0cm\hsize 0cm.76\hsize 0cm.76\hsize 0cm.76\hsize 0cm.76\hsize 0cm.76\hsize 0cm.76\hsize 0cm.74\hsize 0cm.64\hsize 0cm.64\hsize 0cm\hsize 0cm\hsize 0cm \hsize
We shall use the coloring $\chi$ to construct a map $\phi\colon X\to \M$. We will confine the image of $\phi$ to
a single chart of~$\M$; this way we may assume that $M=\R^d$.
We start by mapping the vertices of $X$ to points in general position
\vadjust{\hfill\smash{\raise -34pt\llap{%
\begin{tikzpicture}[round/.style={draw,circle,minimum size=2mm,inner sep=0pt}]
\node[round] (A) at (-1,0) {};
\node[round] (B) at (0,0) {}
   edge [bend right=45] (A);
\node[round] (C) at (1,0) {}
   edge [bend left=45] (B)
   edge [bend right=60] (A);
   \node[align=left] at (0,-0.8)  {\scriptsize $\phi$-image of $X^{\leq (d-1)}$ \\\scriptsize in dimension $d=2$};
\end{tikzpicture}\quad}}}
inside the coordinate hyperplane $\{x_d=0\}$. 
Extending this embedding linearly, we obtain a map $X^{\leq (d-1)}\to \R^d$. 
Then, for each \D{-1}-face $\sigma$ of $X$ we push the interior of $\sigma$
up or down in the last coordinate direction according to the color~$\chi(\sigma)$. The
result of this deformation defines the map $\phi$ on $X^{\leq (d-1)}$.

Finally, we extend the map to the \D-skeleton of $X$ in such a way that if 
$\sigma\in X^{=d}$ and if all the \D{-1}-simplices in $\partial \sigma$
are above (below) the $\{x_d=0\}$ hyperplane, then $\sigma$ is also above (resp.\ below)
that hyperplane.

Let $p\in \R^d$, and suppose that $\phi^{-1}(p)$ intersects all \D-faces of $P_1\ast\dotsb\ast P_{d+1}$
for some $P_i\subset V_i$. Because $\phi$ maps the vertices of $X$ to points in general position and because $\phi$ is
linear on $X^{\leq (d-2)}$, the intersection of $\phi$-images of any $d/2+1$ vertex-disjoint simplices in $X^{\leq (d-2)}$ is empty.
So, $\phi^{-1}(p)$ is contained in the interior of at most $d/2$ many faces of $X^{\leq (d-2)}$.
Let $Q$ be the union of vertex sets of these faces; evidently, $\abs{Q}\leq d^2/2$.

Let $P_i'=P_i\setminus Q$. If $\abs{P_i}\geq 30 (\log n)^{1/(d-1)}$ and $n$ is sufficiently large, then
$\abs{P_i'}\geq 25 (\log n)^{1/(d-1)}$. If $\abs{P_i'}\geq 25 (\log n)^{1/(d-1)}$ for all $i$, then $P_1\ast\dotsb\ast P_{d+1}$ contains monochromatic \D-faces
of each of the two colors. Let $\sigma\in (P_1\ast\dotsb\ast P_{d+1})^{=d}$ be one of these two monochromatic \D-faces.
Since $\sigma$ is monochromatic, the only points of $\phi(\sigma)$ that lie on the hyperplane $\{x_d=0\}$ are those of $\phi(\sigma^{\leq (d-2)})$.
By the definition of $Q$ this implies that $p_d\neq 0$, and hence $\chi(\sigma) = \sign p_d$. Since this holds for monochromatic \D-faces of both colors, we reach a
contradiction. So, $\abs{P_i'} < 25(\log n)^{1/(d-1)}$ for some~$i$, after all. \medskip

\paragraph{Topological intermezzo.} Before tackling the proof of the lower bound from Theorem~\ref{thm:top},
we perform a couple of topological simplifications, and review some definitions and results.

Since the lower bound in Theorem~\ref{thm:top} is about surfaces, we may (and do) assume that $\M$
carries a triangulation. By slight abuse of notation, we use the same symbol $\M$ to denote the triangulation. 
Consider the Boolean algebra generated by all the $\phi$-images of simplices in~$X$.
Let $\rho>0$ be the minimal distance between two disjoint sets in this Boolean algebra. We may assume that the
triangulation of $\M$ is sufficiently fine so that the diameter of each cell is less than~$\rho$. By the simplicial approximation theorem \cite{prasolov}
there is a subdivision $X'$ of $X$ and simplicial map $\phi'\colon X'\to \M$ such that the image
$\phi'(x)$ belongs to the (uniquely defined) simplex of $\M$ whose relative interior contains $\phi(x)$.
If we put $A=\{\sigma\in X : p\in \phi'(\sigma)\}$, then $\bigcap_{\sigma \in A} \phi(\sigma)$ is also non-empty.
Therefore, it suffices to find a topological Pach family for $\phi'$. To avoid extra notation, we continue writing $\phi$
for~$\phi'$, with the added benefit of $\phi$ now being a piecewise linear map.

We say that the map $\phi\colon X\to \M$ is in \emph{general position with respect to $\M$} if
for every simplex $\sigma\in X$ and every simplex $\tau\in \M$, whose dimensions satisfy $\dim \sigma+\dim \tau=d$,
$\phi(\sigma)$ and $\tau$ intersect transversally. By \cite[Ch.~VI]{zeeman}, we may assume (up to small perturbation) that $\phi$ is in general position.

A \emph{$k$-chain in a complex $Y$} (over $\F_2$) is a formal $\F_2$-linear combination of $k$-faces of~$Y$. The $k$\nobreakdash-chains
of $Y$ form an $\F_2$-vector space, which we denote $C_k(Y)$. The \emph{boundary} of a $k$-face $\sigma\in Y$ is the $(k-1)$\nobreakdash-chain
$\partial \sigma$ that is the sum of all $(k-1)$-faces of $\sigma$. The boundary of a general $k$-chain is defined
by linearity, $\partial(\sigma_1+\dotsb+\sigma_m)\eqdef\partial \sigma_1+\dotsb+\partial \sigma_m$.

Given  a map $\phi$ and a pair of chains $a\in C_{d-k}(X)$ and $b\in C_k(\M)$,  we can define the intersection number 
$\phi(a)\cdot b$. If $a$ and $b$ are simplices, let $\phi(a)\cdot b$ be the parity of the number of the intersection points $\phi(a)\cap b$, which is
finite since $\phi$ is in general position. We extend by linearity to the case when $a$ and $b$ are arbitrary chains.
The intersection numbers satisfy
\begin{equation}\label{eq:comm}
  \phi(\partial a)\cdot b=\phi(a)\cdot \partial b;
\end{equation}
see \cite[Eq.~(5) on p.~256]{seifert_threlfall}; for a modern proof of a similar formula for the cap product see
\cite[p.~240 in Section~3.3]{hatcher}.


Finally, we need to state Gromov's overlap theorem \cite{gromov} slightly more precisely. We follow the exposition in 
\cite{dotterrer_kaufman_wagner}. A \emph{$k$-cochain in $X$} is a $\F_2$-linear function on $C_k(X)$. The vector space of
all $k$-cochains is denoted $C^k(X)$. A \emph{coboundary} of a cochain 
$a\in C^k(X)$ is the $(k+1)$-cochain $\delta a$, which is given by $\delta a(\sigma)\eqdef a(\partial \sigma)$. 
Define the norm on $C^k(X)$ by
\begin{equation}\label{eq:norm}
  \norm{a}\eqdef\abs{\supp a}/\abs{X^{=k}},
\end{equation}
where $\supp a$ is the support of~$a$. Up to normalization, this is the Hamming norm.
Since we work with $\F_2$ coefficient, we will occasionally abuse the notation and identify 
cochains with their supports.

A $k$-coboundary is a cochain of the form $\delta a$ for some $(k-1)$-cochain~$a$. 
We say that $a$ \emph{cofills}  $\delta a$. In general, there are many $(k-1)$-cochains $a$ that 
cofill a given coboundary. Cofillings of small norm are of special interest though. We say that $X$ satisfies an \emph{$L$-cofilling 
inequality} in dimension $k$ if for every $k$-coboundary $b$ there is some $a\in C^{k-1}(X)$ such that $\delta a=b$ and $\norm{a}\leq L\norm{b}$.

We say that $X$ is \emph{$\veps$-sparse} if, for every face $\tau\in X$ and every $0\leq k\leq d$, less than $\veps$-fraction of all $k$-simplices 
intersect~$\tau$. 

\begin{theorem}[Gromov's overlap theorem, following \cite{dotterrer_kaufman_wagner}]\label{thm:gromov}
Suppose $X$ is a finite simplicial complex of dimension~$d$, and let $\norm{\cdot}$ be defined by \eqref{eq:norm}. Suppose 
$\phi\colon X\to \M$ is a continuous map into a piecewise linear manifold~$\M$. Suppose that 
\begin{itemize}
\item $X$ satisfies an $L$-cofilling inequality in dimensions $1,\dotsc,d$; and
\item $X$ is $\veps$-sparse; and
\item the cohomology groups of $X$ vanish in dimensions $0,\dotsc,d-1$.
\end{itemize}
Then there is a point $p\in \M$ such that 
\[
  \norm{\{\sigma\in X^{=d} : \phi(\sigma)\cdot p\neq 0 \} } \geq \frac{1}{2(d+1)!L^d}-O(d\veps).
\]
\end{theorem}
(The statement appearing as Theorem~8 of \cite{dotterrer_kaufman_wagner} is slightly different: It is limited to compact connected manifolds,
the condition on the cohomology groups is replaced by a condition on cosystoles, and the conclusion asserts only that
$\phi$-images of a positive fraction of the \D-simplices contain~$p$. The compactness condition can be removed
as explained in Remark 9 of the same paper, the connectedness condition can be removed by restricting to the component containing the image of~$X$, 
the cosystole condition is implied by our cohomology condition, and the stronger conclusion follows from the actual proof of the overlap theorem 
in \cite{dotterrer_kaufman_wagner}.)

For the complex $X=V_1\ast\dotsb\ast V_{d+1}$ of interest in this paper, an $L$-cofilling inequality in dimension~$k$
with $L=\frac{\abs{X^{=k}}}{\abs{X^{=(k-1)}}}\cdot \frac{2^k-1}{n}\leq 2^k\leq 2^d$ was proved in \cite[Proposition 5.8]{dotterrer_kahle}.
An $n$-point discrete space can be thought of as a wedge sum of $n-1$ copies of $\Sph[0]$,
and so $X$ is homotopy equivalent to the wedge sum of $(n-1)^{d+1}$ copies of~$\Sph$, which implies
that the cohomology groups vanish in all dimensions less than~$d$.
Finally, $X$ is clearly $O(1/n)$-sparse.
Therefore, by Theorem~\ref{eq:gromov}, for this $X$ and every continuous map $\phi\colon X\to \M$ we obtain a point $p\in \M$ satisfying
\begin{equation}\label{eq:gromov}
  \norm{\{\sigma\in X^{=d} : \phi(\sigma)\cdot p\neq 0 \} } \geq \frac{1}{(d+1)!2^{d^2+1}}-O(1/n).
\end{equation}

\paragraph{Lower bound in the case $\M\neq \RP^2$.} Let $\widetilde{\M}$ be the universal cover of~$\M$.
Since the complex $X$ is simply-connected, we may lift the map $\phi\colon X\to \M$ to a map $\widetilde{\phi}\colon X\to \widetilde{\M}$ satisfying
$\phi=\prj\circ\, \widetilde{\phi}$, where $\prj\colon \widetilde{\M}\to \M$
is the projection map. Since $\M$ is a manifold, so is $\widetilde{\M}$. Observe that a Pach family for the map $\widetilde{\phi}$
is also a Pach family for the map $\phi$. So, we may replace $\M$ by its universal cover $\widetilde{\M}$ if we wish.

If $\M$ is compact, and is neither $\Sph[2]$ nor $\RP^2$, then it follows
from the classification of compact surfaces that its fundamental group is infinite.
When $\pi_1(\M)$ is infinite, $\widetilde{\M}$ is non-compact because it contains a discrete subset of the form~$\prj^{-1}(p)$. 
In this case, we replace $\M$ by $\widetilde{\M}$ as discussed above.
So, if $M\neq \RP^2$, then we may (and do) assume that $\M$ is non-compact.

Since $X$ is compact but $M$ is not, 
the set $\M\setminus \phi(X)$ is non-empty. By Theorem~\ref{thm:gromov} and the explicit bound \eqref{eq:gromov}, there is a point $p\in \M$ 
such that $\phi(\sigma)\cdot p\neq 0$ for $\frac{1}{192}n^3$ many triangles $\sigma$ of~$X$.
Let $\mathcal{F}$ be the set of these triangles. Let $R\in C_1(\M)$ be a path in the $1$-skeleton of $\M$ from $p$ to a point in $\M\setminus \phi(X)$.
For a triangle $\sigma=x_1x_2x_3\in X$ with $x_i\in V_i$,
let $\sigma_i$ be its edge obtained by removing vertex~$x_i$. Associate to each
$\sigma\in \mathcal{F}$ a vector 
\begin{equation}\label{eq:pivector}
\pi(\sigma)\eqdef \bigl(\phi(\sigma_1)\cdot R,\phi(\sigma_2)\cdot R,\phi(\sigma_3)\cdot R\bigr)\in \F_2^3.
\end{equation}
By the pigeonhole principle, there is $\mathcal{F}'\subset \mathcal{F}$ of size $\abs{\mathcal{F}'}\geq 2^{-3}\abs{\mathcal{F}}$ on which
$\pi(\sigma)$ is constant, say $\pi(\sigma)=\tilde{\pi}$. Because of \eqref{eq:comm}, $\sum_i \tilde{\pi}_i=\phi(\sigma)\cdot \partial R=
\phi(\sigma)\cdot p+0\neq 0$.

Let $H\subset X^{=2}$ be the $3$-partite graph whose edge set consists of all the edges of the triangles $\sigma\in \mathcal{F}'$.
We say that $\tau\in E(H)$ is of \emph{type $i$} if $\tau$ contains no vertex from $V_i$ in $X$ (and hence contains
a vertex from every other $V_j$ in $X$). By the definition of $\mathcal{F}'$ we have $\phi(\tau)\cdot R=\tilde{\pi}_i$
whenever $\tau$ is of type~$i$. Since $\sum \tilde{\pi}_i\neq 0$, it follows that every triangle in $H$
corresponds to a $2$-simplex of $X$ whose $\phi$-image contains~$p$.

The lower bound in Theorem~\ref{thm:top}  now follows from the following result of Nikiforov \cite{nikiforov} (stated here for the
special case that we need). An alternative proof is in \cite[Section~2]{rodl_schacht}.
\begin{theorem}[Main result of \cite{nikiforov}] Let $0<c<\tfrac{1}{2}$. Every $n$-vertex graph with at least $cn^3$ triangles
contains a complete tripartite subgraph with all parts of size $\lfloor c^3 \log n\rfloor$.
\end{theorem}
In our case, the graph has $3n$ vertices and at least $2^{-3}\cdot \tfrac{1}{192}n^3$ triangles, and so we may
take $c=1/(192\cdot 2^3\cdot 3^3)$.

\paragraph{Lower bound in the case $\M=\RP^2$.} In contrast to the case $\M\neq \RP^2$, here we reverse the order of the steps: 
we start by applying Gromov's theorem, and then pass to the universal cover. So, pick $p\in \RP^2$ such 
that $\phi(\sigma)\cdot p\neq 0$ for $\tfrac{1}{192}n^3$ many triangles $\sigma\in X^{=2}$.
Let $\widetilde{\phi}\colon X\to \Sph[2]$ be the lift of $\phi$, and let $\prj\colon \Sph[2]\to \RP^2$ be the projection map. 
Denote by $q_1,q_2$ the preimages of $p$ under~$\prj$. Let $R\in C_1(\Sph[2])$ be a path from
$q_1$ to~$q_2$. 

Because $\widetilde{\phi}$ maps the points of $\phi^{-1}(p)$ to $\{q_1,q_2\}$, 
it is clear that $\phi(\sigma)\cdot p=\widetilde{\phi}(\sigma)\cdot(q_1+q_2)$. Since $\partial R=q_1+q_2$,
it follows that $\phi(\sigma)\cdot p\neq 0$ if and only if $\widetilde{\phi}(\partial \sigma)\cdot R\neq 0$.
We can then introduce a vector $\pi(\sigma)$, defined similarly to \eqref{eq:pivector} but with $\phi$ replaced by~$\widetilde{\phi}$. 
The rest of the argument proceeds unchanged.

\section*{Remarks and open problems}
\begin{enumerate}
\item The condition that the manifold $\M$ has no boundary is not essential. If $\partial \M=B$, we can turn
$\M$ into a manifold without boundary by attaching $B\times [0,1)$ to $\M$ along~$B$.

\item The lower bound in Theorem~\ref{thm:sphere} carries over without any changes to any piecewise linear manifold. However, we believe that the upper bound in Theorem~\ref{thm:top} is sharp for every manifold $\M$ that is not a sphere. 
We are unable to prove this even when $\M$ (or its universal cover) is non-compact. Following the argument in Theorem~\ref{thm:top}, 
this case would follow from the obvious generalization of Nikiforov's result from \cite{nikiforov} to hypergraphs. Alas, this is still an open 
problem, see \cite{rodl_schacht} for partial results.

\item Assuming a suitable hypergraph generalization of Nikiforov's theorem, the lower bound argument in Theorem~\ref{thm:top} generalizes to any non-simply-connected piecewise linear manifold~$M$; one obtains $f_M(M)\geq c_M (\log n)^{1/(d-1)}$ for such manifolds. Indeed, we may assume that $\pi_1(M)$ is finite, as otherwise we may reduce to the non-compact case. Let $\ell$ be the smallest prime divisor of $\abs{\pi_1(M)}$. If $\ell=2$, we consider the universal cover of $M$, and proceed as in the case $M=\RP^2$ using $1$\nobreakdash-chain $R$ which is a sum of $\abs{\pi_1(M)}/2$
paths that form a perfect matching on the preimages of a heavily-covered point. If $\ell>2$, we need to make minor alterations to the proof. 
First, we use the Gromov's overlap theorem with $\F_{\ell}$ coefficients in place of $\F_2$ coefficients; as observed in \cite[Remark 18.2]{dotterrer_kaufman_wagner} Gromov's theorem holds under the condition that the manifold is $\F_{\ell}$-orientable. Since $\ell>2$, the fundamental group of the manifold contains no subgroup of index $2$, and so $M$ is orientable \cite[Proposition 3.25]{hatcher}. To apply the overlap theorem we must verify the cofilling inequality of $X$ with $\F_{\ell}$-coefficients. Fortunately the proof in \cite[Proposition 5.8]{dotterrer_kahle} remains valid for coefficients in any abelian group without any modification (notice that there is a small typo on page 13 of \cite{dotterrer_kahle}, in the formula for $\zeta_v$, the first row on the left hand side should read $\eta_{v_0,v}+d \alpha_{v_0,v}$, instead of $\eta_{v_0,v}+\alpha_{v_0,v}$). So we let $p$ be a heavily-covered (mod $\ell$) point, and we let $q_1,\dotsc,q_N$, where $N=\abs{\pi_1(M)}$, be its preimages in the universal cover. Note that $\ell\mid N$ implies that whenever $\phi(\sigma)\cdot p\neq 0\pmod {\ell}$ then there is a pair of preimages $q_i,q_j$ such that $\widetilde{\phi}(\sigma)\cdot(q_i+q_j)\neq 0 \pmod {\ell}$. Therefore, we may take $R$ to be a path connecting the most popular pair of preimages, and finish as in the proof of Theorem~\ref{thm:top}. 

\item The case of a general compact $\M\neq \Sph[d]$ will require a new argument. For example, we do not know how to
prove the lower bound of $c(\log n)^{1/3}$ for $M=\Sph[2]\times\Sph[2]$ even assuming the generalization of Nikiforov's result
to hypergraphs.

\item We thank Joseph Briggs for valuable feedback on a previous version of this paper.

\end{enumerate}
\bibliographystyle{plain}
\bibliography{pachselection}

\end{document}